\documentstyle[11pt]{article}
\input{amssym.def}
\input epsf

\newtheorem{thm}{Theorem}[section]
\newtheorem{lemma}[thm]{Lemma}
\newtheorem{prop}[thm]{Proposition}

\newtheorem{rmk}[thm]{Remark}

\newenvironment{pf}{\bf Proof:\rm}{\begin{flushright}$\Box$\end{flushright}}
\newcommand{\Qed}{\begin{flushright}$\Box$\end{flushright}}

\def\rank{\mathop{\rm rank}\nolimits}
\def\ker{\mathop{\rm Ker}\nolimits}

\newcommand{\C}{{\Bbb C}}
\newcommand{\D}{{\rm D}}

\newcommand{\R}{{\Bbb R}}

\newcommand{\HH}{{\Bbb H}}

\newcommand{\cB}{{\cal B}}
\newcommand{\cC}{{\cal C}}
\newcommand{\cE}{{\cal E}}

\newcommand{\cH}{{\cal H}}

\newcommand{\cN}{{\cal N}}
\newcommand{\cO}{{\cal O}}

\newcommand{\cU}{{\cal U}}

\newcommand{\rB}{{\rm B}}

\newcommand{\fg}{{\frak g}}

\newcommand{\fk}{{\frak k}}
\newcommand{\bb}{{\backslash\!\backslash}}

\newcommand{\fsl}{\frak{sl}}
\newcommand{\fso}{\frak{so}}
\newcommand{\fsp}{\frak{sp}}

\begin{document}
\bibliographystyle{plain}

\title{Morse Groups in Symmetric Spaces Corresponding to
the Symmetric Group}
\author{Mikhail Grinberg}
\date{December 1, 97}
\maketitle

\section{Introduction}

\subsection{Nilcones in Symmetric Spaces}

Let $\theta : \fg \to \fg$ be an involution of a
complex semisimple Lie algebra, $\fk \subset \fg$
the fixed points of $\theta$, and $V = \fg / \fk$ the
corresponding symmetric space.  The adjoint form
$K$ of $\fk$ naturally acts on $V$.  The orbits and
invariants of this representation were studied by
Kostant and Rallis in [KR].  Let $X = K \bb V$ be
the invariant theory quotient, and $f : V \to X$ be
the quotient map.  The space $X$ is isomorphic to
$\C^r$.  When $\fg = \fg' \oplus \fg'$ and $\theta$ 
acts by interchanging the factors, $K \, | \, V$ is just
the adjoint representation of $G'$.  A detailed study of
the singularities of $f$ in this case leads to the
Springer representations of the Weyl group $W'$ of
$\fg'$ (see [BM], [M], [Sp], [Sl]).  More precisely,
if $P$ is the nearby cycles sheaf of $f$ along the
nilcone $\cN = f^{-1} (0)$, then one obtains
representations of $W'$ on the stalks of $P$.

In [Gr] we studied the singularities of $f$ for an
arbitrary symmetric space (in fact, for a more general
class of objects).  The main result of [Gr] is a
description of the nearby cycles sheaf $P$ in terms of
the Fourier transform.

In this paper, we study the microlocal geometry of $f$
in three particular examples: $V^{I} = \fsl_n$, $V^{II} =
\fsl_n / \fso_n$, and $V^{III} = \fsl_{2n} / \fsp_{2n}$.
Of all the classical symmetric spaces, these are the ones
that have the symmetric group $\Sigma_n$ as their small Weyl
group.  The questions we discuss are equally interesting
for other symmetric spaces, but for the time being,
they remain open outside of the three examples above.
Our main results (Theorems \ref{main1} and \ref{main2})
describe the Morse groups of $P$ with two monodromy
structures (see Section 1.2).  Experts in singularity
theory may also find technical Lemma \ref{curv} to be
of some independent interest.

In the case $V = \fsl_n$, we recover the computation by
Evens and Mirkovi\'c [EM] of the local Euler obstructions
for the nilcone in $\fsl_n$.  We also borrow from [EM]
the idea of exploiting torus symmetry (see our proof of
Lemma \ref{crit_pts}).

In his Ph.D. thesis [Groj], Grojnowski studied a class of
equivariant perverse sheaves on symmetric spaces which is
related to the nearby cycles sheaf $P$ we discuss.  He also
proposed the question of studying the characteristic
varieties of these sheaves.

Discussions with Sam Evens, Ivan Mirkovi\'c, and Ian
Grojnowski have been of great value to me.  I also wish to
thank Robert MacPherson for his guidance and support.

\subsection{Morse Groups of the Nearby Cycles}

We now introduce the geometric setup for studying the
singularities of an algebraic map microlocally
(see [L\^e] for an early appearance of these ideas,
and [Gi], [GM], [KS] for a systematic treatment of
Morse groups).  A more technical discussion will be
given in Section 2.

Let $f : V = \C^d \to X = \C^{n}$ be a dominant
algebraic map, such that $0$ is in the image of
$f$, and is a critical value.  Let $X^{reg} \subset X$
be the set of regular values of $f$ (we do not count
a point $\lambda \in X$ as a regular value if
$f^{-1} (\lambda) = \emptyset$).  Assume that $f$ is a
map without blowing up along the fiber $E = f^{-1} (0)$
(this is a kind of a `well-behavedness' condition; see
[Hi] and Section 2.1 below for a precise definition).
Fix a sufficiently fine stratification of $E$.  Let
$S \subset E$ be a connected stratum.  Associated to the
pair $(f, S)$, there is a non-negative integer $m$
defined as follows.

Fix a point $p \in S$.  Take any smooth function
$g : V \to \R$, such that $p$ is a critical point of the
restriction $g |_S$.  Assume $g$ is generic among all such
functions (more precisely, we need the 2-jet of $g$ at $p$
to be generic).  Fix a small $\lambda \in X^{reg}$ and let
$F_\lambda = f^{-1} (\lambda)$.  Then $m$ is the number of
critical points of $g |_{F_\lambda}$ near $p$.  Note that
if $S$ is open in $E$, and consists of regular points of
$f$, then $m = 1$.

In the language of sheaf theory, the number $m$ is the
multiplicity of the conormal bundle $\Lambda_S = T^*_S X
\subset T^* X$ in the characteristic cycle $SS(P)$ of the
nearby cycles sheaf $P$ of $f$.  It is an important
invariant of $f$.

The multiplicity $m$ is, in fact, the dimension of a vector
space $M_\xi (P)$ which depends on $p$ and $g$ only through
the differential $\xi = d_p g \in \Lambda_S$.  The vector
space $M_\xi (P)$ is called the Morse group of $P$ relative
to $\xi$.  It is defined as follows.

Fix a normal slice $N$ to $S$ through $p$, and small
numbers $0 < \delta \ll \epsilon \ll 1$.  Let
$\rB_\epsilon$ be the $\epsilon$-ball around $p$, and $c$
be the complex codimension of $S$ in $E$.  Then
$$M_\xi (P) = H^c ( N \cap F_\lambda \cap \rB_\epsilon, \,
\{ x \, | \, g (x) = g (p) + \delta \} ).$$
Here, we need to assume that $|\lambda| \ll \delta$,
and that $\xi$ lies in a certain Zariski open subset
$\Lambda^0_S \subset \Lambda_S$, called the set of generic
conormals to $S$.

By construction, there are two commuting monodromy actions
on $M_\xi (P)$.  First, there is an action of the
fundamental group $\pi_1 (\Lambda^0_S)$, coming from the
dependence of the Morse group on $\xi$.  This is called the
microlocal monodromy action.  Second, there is an
action of $\pi_1(X^{reg} \cap \rB_0)$, where $\rB_0$ is a
small ball around the origin in $X$.  This comes from the
choice of $\lambda$; we call it the monodromy in the family
$f$.  The Morse groups $M_\xi (P)$, along with these two
kinds of monodromy, give a great deal of information about
the (perverse) sheaf $P$.  In many specific problems, they
suffice to completely determine the structure of $P$ as an
object in the abelian category of perverse sheaves on $E$.
This, in turn, can be used to analyze other questions about
the singularities of $f$.

\subsection{Statement of Results}

Returning to the situation where $V$ is one of the
symmetric spaces $V^{I, II, III}$, we now identify
the geometric ingredients of Section 1.2.

The fundamental group $\pi_1(X^{reg} \cap \rB_0)$
is the classical braid group $\rB_n$ on $n$ strands.  The
zero fiber $E = \cN$ of the quotient map $f$ is naturally
stratified by $K$-orbits.  These orbits are parametrized
by the partitions of $n$ (in the case $V = V^{II}$, there
is, sometimes, an additional sign parameter).  The
partition corresponding to an orbit $\cO \subset \cN$ is
given by the Jordan normal form of a point in $\cO$
(see Lemma \ref{nilp_orb}).

Fix a partition $\bar n : n = n_1 + \dots + n_k$.  Let
$\cO \subset \cN$ be an orbit corresponding to $\bar n$.
Fix a point $A \in \cO$, and let $\Lambda_A^0 =
\Lambda_\cO^0 \cap T^*_A V$, the set of generic covectors
at $A$.  We can not identify the set $\Lambda^0_A$
explicitly.  Instead, we will work with a certain open
subset $\tilde\Lambda_A \subset \Lambda^0_A$.

Order the numbers $n_i$ so that:
$$\displaylines{{\qquad}
n_1 = \dots = n_{m_1} < n_{m_1 + 1} = \dots =
n_{m_1 + m_2} < \dots \hfill\cr\hfill{}
< n_{m_1 + \dots m_{l-1} + 1} =
\dots = n_{m_1 + \dots + m_l},\qquad\cr}$$
with $m_1+...+m_l=k$.  Let $\rB_{\bar n}$ be the group of
braids on $k$ strands, colored in $l$ colors, with $m_j$
strands of $j$-th color ($j = 1, \dots, l$).  The
following lemma will be proved in Section 3.

\begin{lemma}\label{useful_conormals}
There exists a Zariski open subset $\tilde\Lambda_A
\subset \Lambda^0_A$, and a natural homomorphism
$\rho : \pi_1 (\tilde\Lambda_A) \to \rB_{\bar n}$,
such that $\rho$ is an isomorphism when $V = V^I$ or
$V^{III}$, and a surjection when $V = V^{II}$.
\end{lemma}

Fix a basepoint $\xi \in \tilde\Lambda_A$.  Theorem
\ref{main1} describes the Morse group $M_\xi (P)$ as a
$\rB_n$-module.  Let $\sigma_1, \dots, \sigma_{n-1}$ be
the standard generators of $\rB_n$.  We write $\Sigma_n =
\rB_n / (\sigma_1^2 - 1)$ for the symmetric group on $n$
letters, and $\cH_{-1} (\Sigma_n) = \C [\rB_n] /
(\sigma_1 - 1)^2$ for the Hecke algebra of $\Sigma_n$,
specialized at $q = -1$.  Note that $\cH_{-1} (\Sigma_n)$
has a well defined trivial representation of dimension one,
in which all the $\sigma_i$ act by identity.  We denote
this representation by $1$.

\begin{thm}\label{main1}
{\em \bf (i)}   When $V = V^I$ or $V^{III}$, we have:
$$M_\xi (P) \cong \mbox{\em Ind}_{ \, \Sigma_{n_1} \times
\dots \times \Sigma_{n_k}}^{ \, \Sigma_n} \, 1 \, ,$$
as $\rB_n$-modules.  Here, $\rB_n$ acts on the
right-hand side through the natural homomorphism
$\rB_n \to \Sigma_n$.

{\em \bf (ii)}  When $V = V^{II}$, we have:
$$M_\xi (P) \cong \mbox{\em Ind}_{\, \cH_{-1}
(\Sigma_{n_1}) \times \dots \times \cH_{-1}
(\Sigma_{n_k})} ^{\, \cH_{-1}(\Sigma_n)} \, 1 \, ,$$
as $\rB_n$-modules.  Here, $\rB_n$ acts on the
right-hand side through the natural (semi-group)
homomorphism $\rB_n \to \cH_{-1}(\Sigma_n)$.
\end{thm}

In the case $V = V^{I}$, Theorem \ref{main1} is equivalent
to Theorem 0.2 of [EM].  Theorem \ref{main2} describes
the action of $\pi_1 (\tilde\Lambda_A)$ on $M_\xi (P)$.  Let
$\psi : \rB_{\bar n} \to \Sigma_{n_1} \times \dots \times
\Sigma_{n_k}$ be the natural map.  By part (i) of Theorem
\ref{main1}, in the case $V = V^I$ or $V^{III}$, there is a
natural action $\phi$ of the product $\Sigma_{m_1} \times
... \times \Sigma_{m_l}$ on $M_\xi (P)$, commuting with the
action of $\Sigma_n$.  Let $$\chi: \Sigma_{m_1} \times ...
\times \Sigma_{m_l} \to \{ 1, -1 \}$$
be the character taking a simple transposition in
$\Sigma_{m_j}$ to $(-1)^{n_{m_1+...+m_j}}$.

\begin{thm}\label{main2} {\em \bf (i)}
In the case $V = V^I$ or $V^{III}$, the microlocal
monodromy action of $\pi_1 (\tilde\Lambda_A)$ on
$M_\xi (P)$ is given as $(\phi \otimes \chi) \circ \psi
\circ \rho$.
\end{thm}

For the case $V = V^{II}$, note that there is a natural
map $\zeta : \rB_k \to \rB_n$, obtained by
collecting the $n$ strands into $k$ `ropes,' consisting
of $n_1, \dots, n_k$ strands.  More precisely, if
$\kappa_1, \dots, \kappa_{k-1}$ are a the standard
generators of $\rB_k$, then
$$\zeta (\kappa_i) = \sigma_{m_1 + \dots + m_i} \;
\sigma_{m_1 + \dots + m_i - 1} \; \dots \;
\sigma_{m_1 + \dots + m_{i-1} + 1} \qquad\qquad\qquad$$  
$$\sigma_{m_1 + \dots + m_i + 1} \;
\sigma_{m_1 + \dots + m_i} \; \dots \;
\sigma_{m_1 + \dots + m_{i-1} + 2} \; \dots$$
$$\qquad\qquad\qquad \sigma_{m_1 + \dots + m_{i+1} - 1} \;
\sigma_{m_1 + \dots + m_{i+1} - 2} \; \dots \;
\sigma_{m_1 + \dots + m_i}.$$
The right action of $\rB_n$ on itself descends to an
action $\eta$ of the image $\zeta (\rB_{\bar n})$ on
$M_\xi (P)$ (cf. part (ii) of Theorem \ref{main1}).
Define a homomorphism $\o : \rB_n \to \rB_n$ by $\o :
\sigma_i \mapsto \sigma_i^{-1}$.  Note that
$\o$ preserves the image $\zeta (\rB_{\bar n})$.

\vspace{.1in}

\noindent
{\bf Theorem 1.3 \, (ii)}\label{main2b}
\begin{it}
In the case $V = V^{II}$, the microlocal monodromy
action of $\pi_1 (\tilde\Lambda_A)$ on $M_\xi (P)$ is
given as $\eta \circ \o \circ \zeta \circ \rho$.
\end{it}

\vspace{.1in}

Theorems \ref{main1} and \ref{main2} will be proved in
Section 5.

\section{Geometric Preliminaries}

In this section, we recall the basic definitions
pertaining to Morse groups and nearby cycles, and
prove a technical result (Lemma \ref{curv}) about
the curvature of a general fiber near a point
singularity.

\subsection{Nearby Cycles}

Let $V, X$ be smooth, connected algebraic
varieties over $\C$ with $\dim V = d$ and
$\dim X = r$, and let $f : V \to X$ be a dominant
map.  Write $X^{reg} \subset X$ for the set of
regular values of $f$ (we do not count a point
$\lambda \in X$ as a regular value if $f^{-1}
(\lambda) = \emptyset$).  Let $V^\circ \subset V$
be the preimage $f^{-1} (X^{reg})$; note that it
is a manifold.  Fix a point $x \in f (V) \setminus
X^{reg}$, and let $E = f^{-1} (x)$.  Assume that
$f$ is a map without blowing up along $E$, i.e.,
that there exists a stratification $\cE$ of $E$,
such that for any stratum $S \subset E$, Thom's
${\rm A}_f$ condition holds for the pair
$(S, V^\circ)$.  Recall that the ${\rm A}_f$
condition says that for any sequence of points
$v_i \subset V^\circ$, converging to a limit
$e \in S$, if there exists a limit
$$\Delta = \lim_{i \to \infty} T_{v_i} f^{-1}
(f (v_i)),$$ then $\Delta \supset T_e S$.  This
implies, in particular, that $\dim E = d - r$
(see [Hi] for a detailed discussion of the
${\rm A}_f$ condition).

In this setting, we have a well defined nearby
cycles sheaf $P = P_f$ of the map $f$ along $E$.
It is defined as follows.  Let $U$ be a small
neighborhood of $0$ in $\C$.  Choose an algebraic
arc $\gamma : U \to X$, such that $\gamma (0) = x$,
and $\gamma (\tau) \in X^{reg}$, for $\tau \neq 0$.
We may form the pull-back family $f_\gamma : V_\gamma
\to U$, where $V_\gamma = V \times_{X} U$ and
$f_\gamma$ is the projection onto the second factor.
Set $P_\gamma = \psi_{f_\gamma} \, \C_{V_\gamma}
[d-r]$, the nearby cycles of the functions $f_\gamma$
with constant coefficients (see [KS] for a discussion
of the nearby cycles functor $\psi_g$ for a complex
analytic function $g$).

\vspace{.1in}

\noindent
{\bf Proposition-Definition 2.1}
{\it
$\,$ {\em [Gr, Proposition 2.4]}

{\rm\bf (i)}   The sheaves $P_{f_\gamma}$ for different
$\gamma$ are all isomorphic.  We may therefore omit
the subscript $\gamma$, and call the sheaf
$P_f = P_{f_\gamma}$ the nearby cycles of $f$.  It
is a perverse sheaf on $E$, constructible with
respect to $\cE$.

{\rm\bf (ii)}  The local fundamental group $\pi_1
(X^{reg} \cap \rB_x)$, where $\rB_x \subset X$ is a
small ball around $x$, acts on $P_f$ by monodromy.  We
denote this action by $\mu : \pi_1 (X^{reg} \cap
\rB_x) \to {\rm Aut} (P_f)$.
}

\subsection{Morse Groups}

For a stratum $S \in \cE$, let $\Lambda_S =
T^*_S V \subset T^* V$ be the conormal bundle to $S$.
The conormal variety $\Lambda_\cE \subset T^* V$ to
the stratification $\cE$ is defined by 
$$\Lambda_\cE  = \bigcup_{S \in \cE} \Lambda_S.$$
Let $\Lambda^0_\cE$ be the smooth part of
$\Lambda_\cE$.  Note that $\Lambda^0_\cE \subset T^* V$
is $\C^*$-conic.  For $S \in \cE$ we write
$\Lambda^0_S = \Lambda_S \cap \Lambda^0_\cE$; this
is called the set of generic conormals to $S$.

Any perverse sheaf $R$ on $E$, constructible
with respect to $\cE$, gives rise to a local system
$M (R)$ on $\Lambda^0_\cE$, called the Morse local
system of $R$.  The definition is as follows.  Fix
a stratum $S \in \cE$ and a point $p \in S$.  Let
$\Lambda^0_p = T^*_p V \cap \Lambda^0_S$, and choose
a covector $\xi \in \Lambda^0_p$.  Let $g : V \to \R$
be any smooth function with $g (p) = 0$, and
$d_p g = \xi$.  Fix a normal slice $N \subset V$ to
$S$ through $p$, and let $j : N \cap E \to E$ be the
inclusion.  Choose positive numbers $0 < \delta \ll
\epsilon \ll 1$.  Let $\rB_{p, \epsilon} \subset V$
be the $\epsilon$-ball around $p$ (in some fixed
Hermitian metric), and $c$ be the complex codimension
of $S$ in $E$.  The stalk $M_\xi (R)$ is defined by:
$$M_\xi (R) = \HH^{c-r+d} (E \cap N \cap
\rB_{p, \epsilon}, \, \{ g \geq \delta \}; \, j^* R),$$
where the right-hand side is a relative hypercohomology
group with coefficients in $j^* R$.

Lemma \ref{rank_morse} below identifies the Morse
groups of the nearby cycles sheaf $P = P_f$.  Fix
an algebraic arc $\gamma : U \to X$, and let
$P = P_{f_\gamma}$ (this fixes the `up-to-isomorphism'
ambiguity in the definition of $P$).  For $\tau \in U
\setminus \{0\}$, let $F_{p, \tau} = f_\gamma^{-1}
(\tau) \cap N \cap \rB_{p, \epsilon}$.  The space
$F_{p, \tau}$ is a manifold with boundary.  It is
called a Milnor fiber of $f_\gamma$ at $p$ (see [Mi]).

\begin{lemma}\label{rank_morse}
{\rm\bf (i)}   The Morse group $M_\xi (P)$ may be
identified as follows.  In addition to the choices made
above, pick a number $0 < \tau \ll \delta$.  Then:
$$M_\xi (P) = H^c (F_{p, \tau}, \,
\{ g \geq \delta \}; \, \C),$$ where the right-hand
side is an ordinary relative cohomology group.

{\rm\bf (ii)}  Assume $g|_{\rB_{p, \epsilon}}$ is the
real part of a complex algebraic function $\tilde g : 
\rB_{p, \epsilon} \to \C$.  Then the set $\cC$ of critical
points of $\tilde g|_{F_{p, \tau}}$ is finite.

{\rm\bf (iii)} Write $\cC = \{ C_i \}$, and let $m_i$ be
the multiplicity of the critical point $C_i$ (see {\em
[Mi]} for a definition).  Then we have: $\dim M_\xi (P)
= \sum_i m_i$.
\end{lemma}
\begin{pf}
Part (i) is an immediate consequence of the definitions.
For part (ii), note that Thom's ${\rm A}_f$ condition
and the fact that $\xi$ is in $\Lambda^0_\cE$ imply that
the set $\cC$ does not come near the boundary
$\partial F_{p, \tau}$.  But $\cC$ is an intersection of
an affine variety with a closed ball.  Therefore, we
must have $\dim \cC = 0$.  Part (iii) is an application of
Morse theory.
\end{pf}

\subsection{Curvature of the Nearby Fiber}

Assume now $V \cong \C^d$ is a Hermitian affine space.
Given a covector $\xi \in \Lambda^0_p$, we will denote
by the same letter the corresponding affine functional
$\xi : V \to \C$, with $\xi (p) = 0$.  We then have the
following result about the curvature of the Milnor
fibers $F_{p, \tau}$.

\begin{lemma}\label{curv}
There is a Zarsiki open, dense subset $\Lambda^1_p
\subset \Lambda^0_p$, such that for any compact
$\Delta \subset \Lambda^1_p$ and any $\kappa > 0$, there
exists a number $\tau_0 > 0$, such that for any
$\tau \in \C^*$ with $|\tau| < \tau_0$ and any
$\xi \in \Delta$, every critical point $C$ of
$\xi |_{F_{p, \tau}}$ has the following strong
non-degeneracy property: all eigenvalues of the
Hessian of $\xi |_{F_{p, \tau}}$ at $C$ are greater than
$\kappa$.
\end{lemma}
\begin{pf}
Consider the cotangent bundle $\pi : T^* N \to N$.
Fix a number $0 < \tau_1 \ll \epsilon$.  For each
$\tau \in \C^*$ with $|\tau| < \tau_1$, let
$\Omega_\tau \subset T^* N$ be the conormal bundle
to $F_{p, \tau} \setminus \partial F_{p, \tau}$.  Each
$\Omega_\tau$ is a manifold, obtained by intersecting
an algebraic variety with the cylinder $\pi^{-1}
(N \cap \rB^\circ_{p, \epsilon})$, where 
$\rB^\circ_{p, \epsilon}$ is the interior of
$\rB_{p, \epsilon}$.  Let $\Omega_0 \subset \pi^{-1} 
(\rB^\circ_\epsilon)$ be the limit of the family
$\{ \Omega_\tau \}$ as $\tau \to 0$.  We obtain a
$1$-parameter family $q : \Omega \to \{ |\tau| <
\tau_1 \}$, with $q^{-1} (\tau) = \Omega_\tau$.  The
set $\Omega$ is an intersection of an irreducible
algebraic variety with an open region in $T^* N$.
Note that $\Omega^\circ = \Omega \setminus \Omega_0$
is smooth.

From Lemma \ref{rank_morse}, we see that unless
$M_\xi (P) = 0$ for $\xi \in \Lambda^0_p$, the
fiber $\Omega_0$ contains $T^*_p N$ as an
irreducible component.  When $M_\xi (P) = 0$, the
restriction $\xi |_{F_{p, \tau}}$ has no critical
points, and we have nothing to prove.  Assume now
$M_\xi (P) \neq 0$.  By a result of Hironaka
[Hi, p. 248, Corollary 1] there is an algebraic
stratification of $\Omega_0$ such that for any stratum
$\Sigma \in \Omega$, Thom's ${\rm A}_q$ condition holds
for the pair $(\Sigma, \Omega^\circ)$.  Such a
stratification will contain a stratum $\Sigma_0$ which
is an open subset of $T^*_p N$.  We set:
$$\Lambda^1_p = \{ \xi \in \Lambda^0_p \; | \;
\xi |_{T_p N} \in \Sigma_0 \}.$$
The lemma follows from chasing the meaning of the
${\rm A}_q$ condition.
\end{pf}

\section{Some Linear Algebra}

In this section, we discuss some geometric preliminaries
about the map $f : V \to X$ ($V = V^{I, II, III}$).
The main result here is Proposition \ref{norm_form_conorm},
which gives a normal form for a general conormal vector to
a $K$-orbit in the nilcone.

We regard $V^I = \fsl_n$ as the space of all trace zero
endomorphisms of a vector space $U^I \cong \C^n$.  We regard
$V^{II} = \fsl_n / \fso_n$ as the space of all self-adjoint
trace zero endomorphisms of a vector space $U^{II} \cong \C^n$,
endowed with a non-degenerate quadratic form $\nu$.  Lastly, we
think of $V^{III} = \fsl_{2n} / \fsp_{2n}$ in the following
way.  Let $(U^{III}, \omega)$ be a complex symplectic
$2n$-space.  Then $V^{III}$ is the space of all trace zero
endomorphisms $A$ of $U^{III}$, satisfying $\omega (A u_1, u_2)
= \omega (u_1, A u_2)$, for all $u_1, u_2 \in U^{III}$.  We
will omit the superscripts $^{I, II, III}$, whenever a
statement applies to all three, or when a particular case is
specified.

The map $f$ is described as follows.  Let $V = V^I$ or
$V^{II}$, and $A : U \to U$ be an endomorphism in $V$.
Then the components of the image $f (A) \in X \cong \C^{n-1}$
are given by the elementary symmetric functions in the
eigenvalues $\{ \lambda_1, \dots, \lambda_n \}$ of $A$.
Let now $V = V^{III}$.  Then any endomorphism $A \in V$
has a spectrum with even multiplicities.  Write this
spectrum as $\{ \lambda_1, \lambda_1, \lambda_2, \lambda_2,
\dots, \lambda_n, \lambda_n \}$.  Then the components of the
image $f (A)$ are given by the elementary symmetric functions
in $\{ \lambda_1, \lambda_2, \dots, \lambda_n \}$.

We now describe the $K$-orbits in the nilcone
$\cN = f^{-1} (0) \subset V$.  To each $K$-orbit
$\cO \subset \cN$, we associate a partition $\bar n (\cO)$ of
$n$ as follows.  If $V = V^I$ or $V^{II}$, then $\bar n (\cO)$
is the partition $n = n_1 + \dots + n_k$, where $n_1, \dots, n_k$
are the sizes of the blocks in the Jordan normal form of an
element $A \in \cO$.  When $V = V^{III}$, the Jordan form of
any $A \in \cO$ must have an even number of blocks of each size.
Write $n_1, n_1, n_2, n_2, \dots, n_k, n_k$ for the sizes of
these blocks.  Then $\bar n (\cO)$ is the partition
$n = n_1 + \dots + n_k$.  Given a partition $\bar n : n = n_1 +
\dots + n_k$, we will write $m_1, \dots, m_l$ for the
multiplicities in the collection $\{ n_1, \dots, n_k \}$, as in
Section 1.3.  We will call an orbit $\cO \subset \cN$ (and any
point $A \in \cO$) regular, if $\bar n (\cO)$ is the partition
$n = n$.  Equivalently, $A \in \cN$ is regular whenever
$\rank \, (d_A f) = n - 1$.

\begin{lemma}\label{nilp_orb} {\rm \bf (i)}
When $V = V^I$ or $V^{III}$, the correspondence
$\cO \mapsto \bar n (\cO)$ sets up a bijection between the
$K$-orbits in $\cN$ and the partitions of $n$.

{\rm \bf (ii)}
When $V = V^{II}$, the correspondence $\cO \mapsto \bar
n (\cO)$ gives a surjection from the set of $K$-orbits to
the set of partitions.  If all of the multiplicities
$m_1, \dots, m_l$ for a partition $\bar n$ are even, then
there are two $K$-orbits with $\bar n (\cO) = \bar n$.
Otherwise, there is only one.
\end{lemma}
\Qed

We now discuss the conormal variety $\Lambda$ to the orbit
stratification of $\cN$.  Use the non-degenerate bilinear
form $\rm{tr} (AB)$ on $V$ to identify $V$ and $V^*$, and to
regard $\Lambda$ as a subset of $V \times V \;
(\cong V \times V^* \cong T^* V)$.  The following is a
simple exercise.

\begin{lemma}\label{conorm=commute}
Let $A \in \cN$, and $B \in V$.  Then $(A, B) \in \Lambda$
if and only if $AB = BA$.
\end{lemma}
\Qed

Fix a partition $\bar n : n = n_1 + \dots + n_k$, and an
orbit $\cO \subset \cN$ with $\bar n (\cO) = \bar n$.
Let $\Lambda_\cO \subset \Lambda$ be the conormal bundle to
$\cO$, and $\Lambda^0_\cO \in \Lambda_\cO$ be the set of
generic conormals to $\cO$.  Define $\tilde\Lambda_\cO$ to
be the set of all $(A, B) \in \Lambda_\cO$, such that $B$
has $k$ distinct eigenvalues.

\begin{prop}\label{norm_form_conorm}
{\rm \bf (i)}
The set $\tilde\Lambda_\cO$ is a Zariski open, dense
subset of $\Lambda^0_\cO$.

{\rm \bf (ii)}  Fix a point $\xi = (A, B) \in
\tilde\Lambda_\cO$.  Let $U = U_1 \oplus \dots \oplus U_k$
be the generalized eigenspace decomposition for $B$.
This decomposition is orthogonal with respect to $\nu$
when $V = V^{II}$, and to $\omega$ when $V = V^{III}$.
After a suitable reordering, we have:
$\dim U_i = n_i$ (cases $V = V^I$ or $V^{II}$) or
$\dim U_i = 2n_i$ (case $V = V^{III}$).
The endomorphism $A$ preserves each $U_i$, and the
restriction $A |_{U_i}$ is regular.

{\rm \bf (iii)}  When $V = V^{III}$, there is an
$A$-invariant  decomposition $U_i = U^+_i \oplus U^-_i$,
with $U_i^\pm$ Lagrangian in $U_i$.

{\rm \bf (iv)}  For each $i = 1, \dots, k$, there is
a polynomial $P_i$ of degree $n_i$, such that $B |_{U_i} =
P_i (A |_{U_i})$.
\end{prop}
\begin{pf}
This is an exercise in linear algebra using Lemma
\ref{conorm=commute}.
\end{pf}

\begin{rmk} {\em
The reason the results of this paper do not directly
generalize to other classical symmetric spaces is the
failure of the analog of part (ii) of Proposition
\ref{norm_form_conorm}.}
\end{rmk}

\noindent
{\bf Proof of Lemma \ref{useful_conormals}:}  The set
$\tilde\Lambda_A$ is defined as $\tilde\Lambda_\cO \cap
T_A^* V$.  Let $\Pi_{\bar n}$ be the space of all
unordered $k$-tuples $\tilde u$ of distinct points in $\C$,
colored in $l$ colors, with $m_j$ points of $j$-th color.
Then $\rB_{\bar n} = \pi_1 (\Pi_{\bar n})$.  By Proposition
\ref{norm_form_conorm}, we have a natural map
$r : \tilde\Lambda_A \to \Pi_{\bar n}$ which sends
any pair $(A, B) \in \tilde\Lambda_A$ to the spectrum
of $B$, with each eigenvalue colored according to its
multiplicity.  We define $\rho = r_* :
\pi_1 (\tilde\Lambda_A) \to \rB_{\bar n}$.

It is easy to see that $r$ is a fiber bundle with connected
fiber.  Furthermore, one can check that the fiber is simply
connected when $V = V^I$ or $V^{III}$.  The lemma follows.
\Qed

\section{Identifying the Critical Points}

The main idea of the proof of Theorems \ref{main1} and
\ref{main2} is to pick a generic covector
$\xi = (A, B) \in \tilde\Lambda_\cO$, then compute the
Morse group $M_\xi (P)$ using a normal slice $N$ to
$\cO$ at $A$, which is of a very special form depending
on $B$.

We first describe the construction of $N$ in the case
$V = V^I$.  Fix a covector $\xi = (A, B) \in
\tilde\Lambda_\cO$.  Let $U=\bigoplus_{i=1}^k U_i$ be
the generalized eigenspace decomposition for $B$.
For $1 \leq i, j \leq k$, let $V_{i, j} =
Hom (U_i, U_j)$, so that:
$$V = \bigoplus_{1 \leq i, j \leq k} V_{i, j}.$$ 
Let $T \subset V$ be the parallel translate of the
tangent space $T_A \cO$ through the origin.  One can
check that $T$ splits as:
$$T = \bigoplus_{1 \leq i, j \leq k} T_{i, j},$$
where $T_{i, j} = T \cap V_{i, j}$.  Pick a complement
$\bar{N}_{i, j}$ to $T_{i, j}$ in $V_{i, j}$, and let
$$\bar{N} = \bigoplus_{1 \leq i, j \leq k}
\bar{N}_{i, j}.$$  Then $N$ is the parallel translate of
$\bar{N}$, passing through $A$.  It is easy to check
that $N$ is a normal slice to $\cO$.  We will need
to consider a subset $N_{bd} \subset N$ (``{\em bd}''
stands for block-diagonal), defined as the parallel
translate of $\bar N_{bd} = \bigoplus_{i = 1}^k
\bar N_{i, i}$.

In the case $V = V^{II}$, we construct $N = N^{II}$
using the inclusion $j : V^{II} \to V^I$, which comes
from identifying $U^I$ and $U^{II}$.  The image
$j (V^{II})$ is the anti-fixed points of the involution
$\theta : A \mapsto - A^*$ on $V^I$.  We have $\xi = (A, B)
\in \tilde\Lambda_\cO \subset V^{II} \times V^{II}$.  Apply
the construction described above (in the case $V = V^I$)
for the pair $(j(A), j(B))$, to obtain a normal slice
$N^I \subset V^I$.  Require, in addition, that $N^I$ should
be $\theta$-invariant.  Then set $N^{II} = j^{-1} (N^I)$,
and $N^{II}_{bd} = j^{-1} (N^I_{bd})$.

In the case $V = V^{III}$, we proceed as follows.  Given
a covector $\xi = (A, B) \in \tilde\Lambda_\cO$, let
$U_i, U_i^+, U_i^-$ ($i = 1, \dots, k$) be as in
Proposition \ref{norm_form_conorm}.  For $1 \leq i < j
\leq k$, let:
$$V_{i, j}^{+, +} = (Hom (U_i^+, U_j^+) \oplus
Hom (U_j^-, U_i^-)) \cap V,$$
$$V_{i, j}^{+, -} = (Hom (U_i^+, U_j^-) \oplus
Hom (U_j^+, U_i^-)) \cap V,$$
$$V_{i, j}^{-, +} = (Hom (U_i^-, U_j^+) \oplus
Hom (U_j^-, U_i^+)) \cap V,$$
$$V_{i, j}^{-, -} = (Hom (U_i^-, U_j^-) \oplus
Hom (U_j^+, U_i^+)) \cap V.$$
For $1 \leq i \leq k$, let:
$$V_i^+ = (Hom (U_i^+, U_i^+) \oplus
Hom (U_i^-, U_i^-)) \cap V,$$
$$V_i^{+, -} = Hom (U_i^+, U_i^-) \cap V, \quad
V_i^{-, +} = Hom (U_i^-, U_i^+) \cap V.$$

Let $T \subset V$ be the parallel translate of $T_A \cO$
through zero.  For all the possible values of the
subscript and the superscript, choose a linear complement
$\bar N_*^*$ to $T \cap V_*^*$ in $V_*^*$.
Let $\bar N$ be the direct sum of all the $\bar N_*^*$,
and let $\bar N_{bd} = \oplus_{i = 1}^k \bar N_i^+$.
We define $N$ ($N_{bd}$) to be the parallel translate
of $\bar N$ ($\bar N_{bd}$) through $A$.  This completes
the construction of the normal slice $N$.

Choose a regular value $\tilde\lambda \in X^{reg}$,
corresponding to an $n$-tuple $\{ \lambda_1, \dots
\lambda_n \}$ of distinct eigenvalues
($\sum_{i} \lambda_i = 0$).  For $\tau \in \C$, let
$\tau \cdot \lambda \in X^{reg}$ be given by the
$n$-tuple $\{ \tau \cdot \lambda_1, \dots, \tau \cdot
\lambda_n \}$.  Consider a curve $\gamma: \C \to X$,
defined by $\gamma: \tau \mapsto \tau \cdot \lambda$, and
form the pull-back family $f_\gamma : V_\gamma \to \C$.
Fix a small ball $\rB_{A, \epsilon} \subset V$ around $A$,
and let $F_{A, \tau} = f_\gamma^{-1} (\tau) \cap N \cap
\rB_{A, \epsilon}$, for $0 < |\tau| \ll \epsilon$.
Think of $\xi$ as a linear function on $V$, with
$\xi (A) = 0$

\begin{lemma}\label{crit_pts}
The set of critical points of $\xi |_{F_{A, \tau}}$ is
equal $F_{A, \tau} \cap N_{bd}$.
\end{lemma}
\begin{pf}
By the construction of $N$, any point of the intersection
$F_{A, \tau} \cap N_{bd}$ is critical for
$\xi |_{F_{A, \tau}}$.  To prove the opposite inclusion,
consider first the case $V = V^I$.  Let $H \cong
(\C^*)^{k-1}$ be the torus consisting of all $h \in SL_n$
which act by a scalar on each $U_i$.  The torus $H$ acts
on $V$ as a subgroup of $K = SL_n$, i.e., by conjugation.
By construction, the normal slice $N$ is preserved by this
$H$-action.  So are the covector $\xi$ and the fiber
$f_\gamma^{-1} (\tau)$.

We claim that any critical point $C$ of
$\xi |_{F_{A, \tau}}$ must be fixed by $H$.  To prove this,
assume $C \in F_{A, \tau}$ is a critical point with
$\dim \, H \cdot C \geq 1$.  Then the intersection
$(H \cdot C) \cap F_{A, \tau}$ consists entirely of critical
points of $\xi |_{F_{A, \tau}}$, which contradicts part
(ii) of Lemma \ref{rank_morse}.

It remains to note that the fixed points of $H$
in $F_{A, \tau}$ are precisely the intersection
$F_{A, \tau} \cap N_{bd}$.  This competes the proof for
$V = V^I$.

For $V = V^{II}$, we use the inclusion $j : V^{II}
\to V^I$.  Assume $C \in F_{A, \tau}^{II}$ is a critical
point of $\xi^{II}$.  Then $j (C) \in F_{j(A), \tau}^I$ is
a critical point for the covector $\xi^I : V^I \to \C$,
given by $j(B)$ (this is because the normal slice
$N^I$ is invariant with respect to the involution
$\theta$).  Applying the lemma in the case $V = V^I$,
we find that $j (C)$ is in $N^I_{bd}$ and, therefore,
$C$ is in $N^{II}_{bd}$.

For the case $V = V^{III}$, let $H \cong (\C^*)^k$
be the set of all $h \subset Sp_{2n}$ which act by
a scalar $\kappa_i^\pm$ on each $U_i^\pm$ (note
that we must have $\kappa_i^+ \cdot \kappa_i^- = 1$).
The rest is exactly as in the case $V = V^I$.
\end{pf}

The idea of using torus symmetry to study microlocal
geometry in this setting is due to Evens and
Mirkovi\'c [EM].

The intersection $F_{A, \tau} \cap N_{bd}$ is easy to
describe combinatorially.  Denote by $\cB$ the set of all
maps $\beta : \{ \lambda_1, \dots, \lambda_n \} \to
\{ 1, \dots, k \}$ with $\#\beta^{-1} (i) = n_i$.  We
will sometimes use a shorthand: $\beta (i) =
\beta (\lambda_i)$, for $i \in \{ 1, \dots, n \}$.

\begin{lemma}\label{enum_crit_pts}
For each $\beta \in \cB$, there is a unique endomorphism
$C_\beta \in F_{A, \tau} \cap N_{bd}$, such that $C_\beta$
preserves each vector space $U_i$, and the spectrum of
$C_\beta |_{U_i}$ is equal to $\tau \cdot \beta^{-1} (i)$
(with multiplicities doubled when $V = V^{III}$).  The
assignment $\beta \mapsto C_\beta$ gives a bijection
between $\cB$ and $F_{A, \tau} \cap N_{bd}$.
\end{lemma}
\begin{pf}
This follows from the construction of $N_{bd}$.
\end{pf}

\begin{lemma}\label{morse}
{\bf \em (i)}   Each of the critical points $C_\beta$ is
Morse, i.e., the Hessian of $\xi |_{F_{A, \tau}}$ at
$C_\beta$ is non-degenerate.

{\bf \em (ii)}  We have $\dim M_\xi (P) = \#\cB =
\frac{n!}{n_1! \cdot \dots \cdot n_k!} \,$.
\end{lemma}
\begin{pf}
This is a general position argument.  However, we need
to be careful, because the pair $(\xi, N)$ is far from
generic in our construction.

Fix an affine normal slice $N_0$ to $\cO$ at $A$.
Apply Lemma \ref{curv} to the family $f_\gamma$
and the normal slice $N_0$ to obtain an open set
$\Lambda^1_A \subset \Lambda^0_A$.

\vspace{.1in}

\noindent
{\bf Claim:}  Assume $\xi \in \Lambda^1_A$.  Then,
for $0 < |\tau| \ll \epsilon$, all the critical
points of $\xi |_{F_{A, \tau}}$ are Morse.

\vspace{.1in}

To prove the claim, we use the $K$-action.  Let
$\fk_A \subset \fk = \mbox{Lie} (K)$ be the stabilizer
of $A$.  Choose a complement $\fk_A^\perp$ to $\fk_A$ in
$\fk$.  Assuming the number $\epsilon > 0$ is sufficiently
small, there exists a closed neighborhood $\cU \subset N_0$
of $A$, and a unique diffeomorphism $\psi : \cU \to N \cap
\rB_\epsilon$, such that $y = \psi (x)$ if and only if
$y = exp(t) \, x$, for some $t \in \fk^\perp_A$.  Note
that $f \circ \psi = f$, so that $F_{A, \tau} =
\psi (f_\gamma (\tau) \cap \cU)$.

Choose a compact neighborhood $\Delta \subset \Lambda^1_A$
of $\xi$.  Note that the tangent spaces $T_A N_0$ and
$T_A N$ are naturally identified with the quotient
$T_A V / T_A \cO$, and therefore with each other.  Under
this identification, the differential $d_A \psi$ is the
identity.  It follows that we can choose a $\tau_1 > 0$,
such that for any $\tau \in \C^*$ with $|\tau| < \tau_1$,
and any critical point $C$ of $\xi |_{F_{A, \tau}}$,
we have: $d^*_{\psi^{-1} (C)} \psi \, (\xi) \in \Delta$
(we use the identifications $T^*_C N \cong
T^*_{\psi^{-1} (C)} N_0 \cong (T_A V / T_A \cO)^*$).
Next, choose a number $\kappa > 0$ which is large compared
to the second derivatives of $\psi$.  With these choices,
we may use Lemma \ref{curv} to select a $\tau_0 > 0$.  The
claim for $|\tau| < \min\{\tau_0, \tau_1\}$ then follows
from Lemma \ref{curv}.

The claim, together with Lemmas \ref{rank_morse},
\ref{crit_pts}, and \ref{enum_crit_pts}, implies that
$\dim M_\xi (P) = \#\cB$, if $\xi \in \Lambda^1_A$.
However, $\dim M_\xi (P)$ is independent of $\xi$, for
$\xi \in \Lambda^0_A$.  Applying Lemma \ref{rank_morse}
in the other direction, we conclude that for any
$\xi \in \tilde\Lambda_A$, all the critical points
$C_\beta$ are Morse.
\end{pf}

\section{Picard-Lefschetz Theory}

Given Lemmas \ref{crit_pts}, \ref{enum_crit_pts}, and
\ref{morse}, the proofs of Theorems \ref{main1} and
\ref{main2} are obtained as applications of
Picard-Lefschetz theory.

Continuing with the situation of Section 4, assume
that the numbers $\{ \lambda_1, \dots, \lambda_n \}$
satisfy $\lambda_1 < \dots < \lambda_n$.  Assume
also that the covector $\xi = (A, B)$ is chosen so
that the endomorphism $B$ is semisimple, with
eigenvalues $u_1 < \dots < u_k$, and that
$B |_{U_i} = u_i$.

Let $\Pi_n$ denote the set of all unordered $n$-tuples
of distinct points in $\C$.  We think of  $\rB_n$ as
the fundamental group $\pi_1 (\Pi_n, \tilde\lambda)$,
where $\tilde\lambda = \{ \lambda_i \}$.  Then
$\Sigma_n$ is the permutations of the $\{ \lambda_i \}$,
and it acts on the index set $\cB$ by
$w : \beta \mapsto w \beta = \beta \circ w^{-1}$.  Let
$\beta_0 \in \cB$ be the unique map with $\beta_0
(\lambda_i) \leq \beta_0 (\lambda_j)$, for all
$1 \leq i < j \leq n$.  Every $\beta \in \cB$ is of the
form $w \beta_0$, for some $w \in \Sigma_n$.

Fix small numbers $0 < \tau \ll \delta \ll \epsilon \ll 1$,
define $F_{A, \tau} = f_\gamma^{-1} (\tau) \cap N
\cap \rB_{A, \epsilon}$, and recall (cf. Lemma
\ref{rank_morse}) that
$$M_\xi (P) = H^c (F_{A, \tau}, \,
\{ \xi \geq \delta \}; \, \C),$$
where $c$ is the complex codimension of $\cO$ in $\cN$.

We now recall a standard Picard-Lefschetz construction
of classes in $M_\xi (P)$ (see [AGV] for a detailed
discussion of Picard-Lefschetz theory).  Let $C = C_\beta$ be
one of the critical points of $\xi |_{F_{A, \tau}}$ described
in Lemma \ref{enum_crit_pts}.  Note that
$\xi (C) = \tau \cdot \sum_{i=1}^n \sigma_i \cdot
u_{\beta(i)}$.  Fix a smooth path $\alpha : [0, 1] \to \C$
such that:

(i)   $\;\; \alpha (0) = \xi (C)$, and $\alpha (1) =
      \delta$;

(ii)  $\; |\alpha (t)| \leq \delta$, for any $t \in [0, 1]$;

(iii) $\alpha (t) \neq \xi (C_{\beta'})$, for any
      $t > 0$, $\beta' \in \cB$;

(iv)  $\alpha (t_1) \neq \alpha (t_2)$, for $t_1 \neq t_2$;

(v)   $\; \alpha' (t) \neq 0$, for $t \in [0, 1]$.

\noindent
Let $\cH : T_C F_{A, \tau} \rightarrow \C$ be
the Hessian of $\xi \, |_{F_{A, \tau}}$ at $C$, and
let $T_C \, [\alpha] \subset T_C F_{A, \tau}$ be the
positive eigenspace of the (non-degenerate) real quadratic
form $$\mbox{Re} \, (\cH / \alpha'(0)) :
T_C F_{A, \tau} \to \R.$$  Note that $\dim_\R \, T_C \,
[\alpha] = \dim_\C \, F_{A, \tau} = c$.  Fix an orientation
$\cO$ of $T_C \, [\alpha]$.  The triple $(C, \alpha, \cO)$
defines a homology class
$$[C, \alpha, \cO] \in H_c (F_{A, \tau}, \,
\{ \xi \geq \delta \}; \, \C) = M_\xi (P)^*.$$
Namely, the class $[C, \alpha, \cO]$ is represented by an
embedded $c$-disc
$$\kappa : (\D^c, \, \partial \, \D^c) \to
(F_{A, \tau}, \, \{ \xi \geq \delta \}),$$
such that the image of $\kappa$ projects onto the image of
$\alpha$ and is tangent to $T_C \, [\alpha]$ at $C$.  The
sign of $[C, \alpha, \cO]$ is given by the orientation
$\cO$.  It is a standard fact that
$[C, \alpha, \cO] \neq 0$.

Note that a choice of $\sqrt{-1}$ gives an isomorphism
$M_\xi (P) \cong M_{-\xi} (P)$.  On the other hand, the
intersection pairing on $F_{A, \tau}$ induces a perfect
duality between $M_\xi (P)$ and $ M_{-\xi} (P)$.  Therefore,
the Morse group $M_\xi (P)$ is canonically its own dual,
and we may regard the Picard-Lefschetz class
$[C, \alpha, \cO]$ as an element of $M_\xi (P)$.

Before we begin the proofs of Theorems \ref{main1} and
\ref{main2}, we need to recall a result of [Gr].

\begin{thm}\label{quote} {\em [Gr, Theorems 3.1, 6.4]}

{\rm\bf (i)}  In the case $V = V^I$ or $V^{III}$, the
monodromy action $\mu : \rB_n \to \mbox{Aut} (P)$ factors
through $\Sigma_n$.

{\rm\bf (ii)}  In the case $V = V^{II}$, the monodromy
action $\mu : \rB_n \to \mbox{Aut} (P)$ factors through
$\cH_{-1} (\Sigma_n)$.
\end{thm}

In the case $V = V^{I}$, Theorem \ref{quote} goes back to
the work of Slodowy in [Sl], and appears in [M] in its
present form.

\vspace{.1in}

\noindent
{\bf Proof of Theorem \ref{main1}:}
Begin with part (i) of the theorem.  By part (i) of
Theorem \ref{quote}, the monodromy action $\mu_* :
\rB_n \to \rm{Aut} (M_\xi (P))$ factors through an action
$\tilde\mu : \Sigma_n \to {\rm Aut} (M_\xi (P))$.  Write
$C_0 = C_{\beta_0}$.  Let $\alpha_0 : [0, 1] \to \C$ be the
straight line path connecting $\xi (C_0)$ to $\delta$.  Pick
an orientation $\cO_0$ of $T_{C_0} \, [\alpha_0]$.  Let
$e_0 = [C_0, \alpha_0, \cO_0] \in M_\xi (P)$.  Consider
the product $\Sigma_{n_1} \times \dots \times \Sigma_{n_k}$
as a subgroup of $\Sigma_n$ in the obvious way.

\vspace{.1in}

\noindent
{\bf Claim 1:}  The image of $e_0$ under the action
$\tilde\mu$ spans $M_\xi (P)$ as a vector space.

\vspace{.1in}

\noindent
{\bf Claim 2:}  For $w \in \Sigma_{n_1} \times
\dots \times \Sigma_{n_k} \subset \Sigma_n$, we have:
$\tilde\mu (w) \, e_0 = e_0$.

\vspace{.1in}

Part (i) of the theorem follows immediately from the
two claims and the computation of $\dim \, M_\xi (P)$
in Lemma \ref{morse}.  The proofs of both claims are
standard applications of Picard-Lefschetz theory.  Begin
with Claim 1.  For each $\beta \in \cB$, pick an
element $w_\beta \in \Sigma_n$ with $w_\beta \,
\beta_0 = \beta$.  Then $\tilde\mu (w_\beta) \, e_0$
is of the form $[C_\beta, \alpha_\beta, \cO_\beta]$,
where $\alpha_\beta$ is some path connecting
$\xi (C_\beta)$ to $\delta$, and $\cO_\beta$ is an
orientation of $T_{C_\beta} \, [\alpha_\beta]$.  The
set $\{ \tilde\mu (w_\beta) \, e_0 \}_{\beta \in \cB}$
contains one Picard-Lefschetz class for each critical
point of $\xi |_{F_{A, \tau}}$.  Therefore, it is a
basis of $M_\xi (P)^*$.

To prove Claim 2, let $s$ be a simple reflection in
$\Sigma_{n_1} \times \dots \times \Sigma_{n_k}$.
Note that $s \, \beta_0 = \beta_0$.  By tracing what
happens to the critical values $\xi (C_\beta)$ as we
permute the two eigenvalues corresponding to $s$, it is
not hard to see that the path $\alpha_0$ does not
change, and we have: $\tilde\mu (s) \, e_0 = \pm e_0$.
A further argument with the Hessian $\cH_0 : T_{C_0}
F_{A, \tau} \to \C$ of $\xi |_{F_{A, \tau}}$ at $C_0$
shows that, in fact, $\tilde\mu (s) \, e_0 = e_0$.

The analysis of $\cH_0$ is based on the fact that the
decomposition of $T_{C_0} F_{A, \tau}$ by the
intersections with the subspaces $\bar N_{i,j} \subset
T_{C_0} V \cong V$ used in the construction of $N$ (see
Section 4) has the following properties:

(i)  $T_{C_0} F_{A, \tau} \cap \bar N_{i,i} = 0$, and

(ii) $T_{C_0} F_{A, \tau} \cap \bar N_{i,j}$ is
orthogonal to $T_{C_0} F_{A, \tau} \cap \bar N_{l,m}$
with respect to $\cH_0$, unless $i = m$ and $j = l$.

\noindent
We omit the details of the analysis of $\cH_0$.

Part (ii) of the theorem is proved similarly.
\Qed

\vspace{.1in}

\noindent
{\bf Proof of Theorem \ref{main2}:}
This is similar to the proof of Theorem \ref{main1}.
Continuing with the notation of that proof, we give an
outline of the argument.

{\bf Step 1.}  Denote the microlocal monodromy action of
$\pi_1 (\tilde\Lambda_A)$ on $M_\xi (P)$ by $h$.  Since $h$
commutes with the monodromy action $\mu_*$, and $e_0$
generates $M_\xi (P)$ under $\mu_*$, it is enough to verify
the claim of the theorem about $h (\upsilon) \, e_0$, for
any $\upsilon \in \pi_1 (\tilde\Lambda_A)$.

{\bf Step 2.}  We show that $h$ factors through the
homomorphism $\rho$ of Lemma \ref{useful_conormals}.
In the case $V = V^I$ or $V^{III}$, this is clear since
$\rho$ is an isomorphism.  Consider now the case
$V = V^{II}$.  Choose an element $\upsilon \in \ker
\rho$.  Recall that the homomorphism $\rho :
\pi_1(\tilde\Lambda_A) \to \rB_{\bar n}$ is defined as 
the push-forward $r_*$ by the map $r : \tilde\Lambda_A \to
\Pi_{\bar n}$ which takes a covector $(A, B)$ to the
spectrum of $B$.  Therefore, $\upsilon$ may be represented
by a loop lying inside a fiber of $r$.  It follows that
$h (\upsilon)$ fixes any Picard-Lefschetz class up to sign.
A verification with the Hessian $\cH_0$ shows that, in fact,
$h (\upsilon) \, e_0 = e_0$.  We may therefore write
$h = \tilde h \circ \rho$.

{\bf Step 3.}  We now specify a set of generators for
the group $\rB_{\bar n}$.  Let $\tilde u = r (\xi)
= \{ u_1, \dots, u_k \} \in \Pi_{\bar n}$.  We identify
$\rB_{\bar n}$ with the fundamental group
$\pi_1 (\Pi_{\bar n}, \tilde u)$.  Note that $\rB_{\bar n}$
is naturally a subgroup of $\rB_k$.  Let $\kappa_1, \dots,
\kappa_{k-1}$ be the standard generators for $\rB_k$.
For each $1 \leq i < j \leq l$, define an element
$\varsigma_{i,j} \in \rB_{\bar n}$ by
$$\varsigma_{i, i+1} = \kappa^2_{m_1 + \dots + m_i},
\;\; \mbox{and}$$
$$\varsigma_{i,j} = \kappa^{-1}_{m_1 + \dots + m_{j-1}} \;
\kappa^{-1}_{m_1 + \dots + m_{j-1} - 1} \; \dots \;
\kappa^{-1}_{m_1 + \dots + m_i + 1} \;
\kappa^2_{m_1 + \dots + m_i} \qquad\qquad$$
$$\qquad\qquad\qquad\qquad\qquad\qquad
\kappa_{m_1 + \dots + m_i + 1} \;
\kappa_{m_1 + \dots + m_i + 2} \; \dots \;
\kappa_{m_1 + \dots + m_{j-1}},$$
when $j - i > 1$.  Note that $\varsigma_{i,j}$ is just a
braid moving the point $u_{m_1 + \dots + m_{j-1} +1}$
once around $u_{m_1 + \dots + m_i}$.  It is not hard to
check that the $\{\varsigma_{i,j}\}$, together with
the $\{ \kappa_i \; | \; \kappa_i \in \rB_{\bar n} \}$
give a set of generators for $\rB_{\bar n}$.

{\bf Step 4.}  Let $\upsilon$ be one of the generators of
$\rB_{\bar n}$ constructed in Step 3.  By tracing what
happens to the path $\alpha_0$ used in the definition of
$e_0$, as we vary $\tilde u \in \Pi_{\bar n}$ along a path
representing $\upsilon$, we can check that
$$\tilde h (\upsilon) \, e_0 = \pm \,
\mu_* \circ \o \circ \zeta \, (\upsilon^{-1}) \, e_0.$$

{\bf Step 5.}  A verification with the Hessian $\cH_0$ shows
that the sign in the above equation is a plus, except when
$V = V^I$ or $V^{III}$, and $\upsilon \in
\{ \kappa_i \; | \; \kappa_i \in \rB_{\bar n} \}$.

{\bf Step 6.}  Recall that under the identification of
$M_\xi (P)$ with a quotient of $\C [\rB_n]$ in the proof
of Theorem \ref{main1}, the element $e_0$ is the image of
$1 \in \C [\rB_n]$.  Therefore, Steps 4 and 5 verify the
claim of the theorem about the action of
$\tilde h (\upsilon)$ on $e_0$.
\Qed

\vspace{.1in}

\noindent
Department of Mathematics, MIT, 77 Massachusetts Ave.,
Cambridge, MA 02139

\noindent
{\it grinberg@math.mit.edu}

\end{document}